\documentclass[10pt,a4paper]{article}
\usepackage{typearea}
\usepackage{amsfonts}
\usepackage{amssymb}
\usepackage{amsthm}
\usepackage{amscd}
\usepackage[centertags]{amsmath}
\usepackage{eucal}
\usepackage{epsfig}
\usepackage[matrix,arrow,curve]{xy}
\CompileMatrices
\typearea{12}
\author{Andreas Klein}
\title{Hamiltonian fixed points, symplectic spinors and Frobenius structures}

\newtheorem{theorem}{Theorem}[section]

\newtheorem{folg}[theorem]{Corollary}

\newtheorem{conj}[theorem]{Conjecture}

\begin{document}
\maketitle

\begin{abstract} This article announces a series of articles aiming at introducing the concept of symplectic spinors into symplectic topology resp. the concept of Frobenius structures. We will give lower bounds for the number of fixed points of a Hamiltonian diffeomorphism on the cotangent bundle over a compact manifold $M$ by defining a certain $C^*$-valued function on $T^*\tilde M$, where $\tilde M$ is a certain 'complexification' of $M$, whose critical points are closely related to the fixed points of the Hamiltonian diffeomorphism $\Phi$ in question. This function, defined via embedding $\tilde M$ into $\mathbb{R}^m$ for an appopriate $m$ and the use of symplectic spinors, is essentially determined by associating to each point of $T^*\tilde M$ the value of a certain spinor-matrix coefficient of specific elements of the Heisenberg group which are determined by $\Phi$. We will discuss an approach for the case of the torus $M$ which does not require embeddings. Here, the matrix coefficients in question coincide with a certain theta function associated to the Hamiltonian diffeomorphism. We will discuss how to define spectral invariants in the sense of Viterbo and Oh by lifting the above function to a real-valued function on an appropriate cyclic covering of $T^*\tilde M$ and using minimax-methods for 'half-infinite' chains. Furthermore we will define a 'Frobenius structure' on $T^*\tilde M$ by letting elements of $T(T^*\tilde M)$ act on the fibres of a line bundle $E$ on $T^*\tilde M$ spanned by 'coherent states' closely related to the above spinor-matrix coefficient. The spectral Lagrangian in $T^*(T^*\tilde M)$ associated to this Frobenius structure intersects the zero-section $T^*\tilde M$ exactly at the critical points of the function described beforehand.

\end{abstract}

\section{Research announcement}\label{intro}

This announces a series of articles (\cite{kleinlag}, \cite{kleinham1}, \cite{kleinham2}) which aim to introduce the concept of symplectic spinors (Kostant \cite{kost}) into symplectic topology on one hand and the field of 'Frobenius structures' as introduced by Dubrovin (\cite{dubrovin}) on the other hand. Note that neither the former nor the latter relation is completely new in the mathematical literature, as can be read off for instance from the occurence of symplectic spinors in the literature concerning the Maslov index, semiclassical approximation and geometric quantization (cf. Guillemin, Leray, Crumeyrolle \cite{cru1}, \cite{guillemin}, \cite{leray}) on one hand and the introduction of the 'Geometric Weil representation' by Deligne (letter to Kazhdan, 1982 \cite{deligne}) on the other hand. The latter was reinforced in contemporary discourse in the realms of the Langlands program (cf. Lafforgue and Lysenko \cite{lafforgue}) resp. the 'mirror-symmetry'-conjecture first introduced by Kontsevich into mathematics. However, as far as the author knows, there has been no systematic treatment yet to explore the possible role of the notion of symplectic spinors and the Weil representation in 'modern symplectic topology', which can be traced back to pseudoholomorphic curves introduced by Gromov and the advent of infinite dimensional variational methods as introduced by Floer. In between both, one can consider the finite dimensional variational methods of Viterbo (\cite{viterbo}) and their relation to symplectic capacities as introduced by Hofer (\cite{hofer}) and exactly this will be the starting point of this series of papers. The main observation linking symplectic spinors to symplectic topology on one hand and 'Frobenius structures' on the other hand is the existence of a construction which links Lagrangian submanifolds of the cotangent bundle $T^*M$ of a compact Riemanian manifold $M$, intersecting each cotangent fibre transversally, at least outside of their 'caustic' to sums of complex lines, viewed as subbundles in the symplectic spinor bundle, that is we have a correspondence:
\[
{\rm (unramified)\ Lagrangian\ submanifolds\ of}\ T^*M \quad \leftrightarrow\quad {\rm direct\ sums}\ \bigoplus_i (\mathcal{L}_i\rightarrow M)
\]
where $\mathcal{L}_i, i=1,\dots, k$ are a certain set of complex line-subbundles of the symplectic spinor bundle $\mathcal{Q}$ on $T^*M$, pulled back to $M$, $i^*\mathcal{Q}$, where $i:M\hookrightarrow T^*M$ is the inclusion of the zero section. Recall that the symplectic spinor bundle $\mathcal{Q}$ over the symplectic manifold $(T^*M^{n},\omega)$ is the bundle associated to a certain connected $2$-fold covering of the principal bundle of symplectic frames, called a metaplectic structure, by the Shale-Weil-representation of the connected $2$-fold cover of the symplectic group acting as intertwining operators for the Schroedinger representation $\rho$ of the Heisenberg group $H_n$ on $L^2(\mathbb{R}^n)$. Metaplectic structures exist under relatively mild conditions on $M$, that is if $c_1(T^*M)=0\ {\rm mod}\ 2$. Note that each branch of the Lagrangian submanifold $\pi:L\subset T^*M\rightarrow M$ covering $M$ gives over any $x \in M$ rise to an element $\psi_{i,x} \in i^*\mathcal{Q}_x\simeq L^2(\mathbb{R}^n)$ by setting
\[
\psi_{i,\lambda,x}(u)=\rho((0,p_i),\lambda)f(u), \quad ((0,p),\lambda) \in H_n, \ u\in \mathbb{R}^n.
\]
Here, $p_i \in \mathbb{R}^n$ locally parametrizes the $i$-th branch of $L$, $\lambda \in \mathbb{R}$ (arbitrary at this point) and $f \in L^2(\mathbb{R}^n)$ is the Gaussian, we identify $H_n=\mathbb{R}^{2n} \times \mathbb{R}$. The set $\psi_{i,\lambda, x}, \ x \in M$ defines a smooth complex line bundle $\mathcal{L}_i$ (outside of ramification points) over $M$ since $i^*\mathcal{Q}_x$ allows a reduction to the structure group $O(n)$ (or its two-fold covering) and $\rho$ acts equivariantly w.r.t. to the Shale-Weil-representation. By construction, $k$ equals the local number of branches of $L$. Note that physically, the vectors $\psi_{i,\lambda,x}$ correspond exactly to 'coherent states' of the quantum mechanical Harmonic oscillator. The above correspondence will be called a {\it symplectic Fourier Mukai transformation}. In this and the second paper in this series, we will mostly assume that $\pi$ is of constant non-zero degree (hence surjective) and the set of caustic points ${\rm ker}\ d\pi\cap TL\neq\{0\}$ is empty. Under this hypothesis, each branch of the above non-ramified Lagrangian furthermore corresponds to a summand of a certain $\mathbb{C}$-valued function on $E$, where $E:= \bigoplus_i^k \mathcal{L}_i\rightarrow M$, namely we pair the above $\psi_{i,\lambda_i, x} \in i^*\mathcal{Q}_x$ over each point $x \in M$ with certain 'elementary vectors'(cf. \cite{mumford}). Let us assume each fibre $E_x$ carries a lattice $\Gamma_x$ being compatible with $L\cap E_x$ in the sense that $L=p^{-1}(\tilde L)$ for a Lagrangian $\tilde L$ in the torus bundle $p:E\rightarrow E/\Gamma$. Then, by duality, the structure group of $i^*\mathcal{Q}$ is reducible to $O(n)\cap Sp(2n, \mathbb{Z})$. In this situation, the canonical pairing in $i^*\mathcal{Q}$ of the  $\psi_{i,\lambda_i, x}$ with another (the globally defined) distinguished vector $e_\mathbb{Z} \in i^*\mathcal{Q}_x$, which can be considered as a sum of delta distributions centered on the integer points of $\mathbb{R}^n$, defines over each point of $E$ a sum of matrix elements, that is a mapping
\[
\Theta: E\rightarrow \mathbb{C},\quad (x,c) \mapsto \sum_i^k<\psi_{i,\lambda_i, x}, e_\mathbb{Z}>(c)
\]
where we extend over each fibre $\mathcal{L}_{i,x}$ by multiplying the argument of $\rho$ acting on $f$ as well as the argument of $e_\mathbb{Z}$ by an affine-linear polynomial in $c$ ($c=(c_i)_{i=1}^k$ is the complex coordinate on the fibres of $E$, for details see \cite{kleinlag}). In case of {\it exact} $L$, that is, the canonical one-form $\alpha$ on $T^*M$ is exact on $L$, we will fix the above $\lambda_i$ by being the integral of the Poincare-Cartan-form $\alpha_H=\alpha-H_tdt$ along rays emanating from $x$ to the $i$-th branch of $L$, where $H_t$ is defined so that its Hamiltonian flow generates these rays. Choosing an appropriate basis for $i^*T(T^*M)$, each summand of this function, evaluated at $x \in M$, considering $M$ as the zero-section of $E$, can be interpreted as a value of a certain (sum of) theta functions, that is of functions of the form
\[
\theta(z, \Omega)=\sum_{k \in \mathbb{Z}^n}e^{\pi i (k, \Omega k)+2 \pi i(k, z)+i\lambda},
\]
where $\Omega$ is an element of the Siegel upper half space (a symmetric complex $n\times n$-matrix $\Omega$ whose imaginary part is positive definite) and $(\cdot,\cdot)$ denotes the standard sesquilinar form on $\mathbb{C}^n$. Note that in the case the above torus-bundle structure is absent, we will use different distinguished vectors of $i^*\mathcal{Q}_x$ to define $\Theta$, one choice is to replace $e_\mathbb{Z}$ by the Gaussian $f$. The above choice $e_\mathbb{Z}$ in the presence of a transversal Lagrangian $L$ and a compatible lattice $\Gamma$ will be considered as the most fundamental for reasons that will hopefully become clearer in the course of this article and its followers. To summarize the above philosophically, we want to stress that using these constructions, there is a local correspondence between Lagrangian submanifolds and (special values of) theta functions on one hand and complex line bundles over $M$ on the other hand, as long as the latter are spanned by 'coherent states'. For this terminology, see Perelmov (\cite{perelmov}). If $L$ is furthermore exact, then choosing the data as above, $\Theta$, outside of an eventual zero set $S$ (to be interpreted as some sort of theta divisor) defines a generating function $\Theta:E\setminus S\rightarrow \mathbb{C}^*$ for $L$ (generalizing Viterbo's construction) that reproduces $L$ up to scaling by a non-zero complex number in the fibres of $T^*_{\mathbb{C}}M$ (the 'logarithmic derivative') and, lifted to a suitable cyclic covering $\tilde E$ (associated for instance to $\Theta_*:\pi_1(E\setminus S)\rightarrow \pi_1(S^1)$), allows to define spectral invariants in a very similar way, using the Morse theory for Novikov one forms developed by Novikov, Farber, Ranicki and others. The critical points of $\Theta$ then correspond to the intersection points of $L$ with the zero section. Note that $\tilde E$ is a vector bundle over a (non-compact) cyclic covering $\tilde M$, of $M$.\\ 
Finally, since the vectors $\psi_{i,\lambda_i}$ define a non-vanishing section of $E=\bigoplus_i^k  \mathcal{L}_i$ on $M$, symplectic Clifford multipliction on $T^*M$ allows us to define a Frobenius multiplication $\star$ in the sense of Dubrovin \cite{dubrovin} for tangent vectors on $M$, that is for $v \in TM$ we set
\[
\star \in H^0(T^*M\otimes End(E)),\quad v\star \psi_i:=i(v-iJv)\cdot\psi_i,
\]
where $\cdot$ denotes symplectic Clifford multiplication over $T^*M$ and $J$ denotes a compatible nearly complex structure on $T(T^*M)$. As it turns out, the $\psi_i$ diagonalize $\star$ and its eigenvalues ($\star$ is semisimple, which is a consequence of our assumption of $L$ being non-ramified), considered as elements of $\Gamma(\Lambda^1(T^*M))$, are precisely the branches of the above Lagrangian submanifold $L$, that is, we recover $L$ as the spectral Lagrangian of $\star$. As a variety, this Lagrangian thus identifies with 
\[
L\simeq {\rm Spec}(\frac{{\rm Sym}(TM)}{\mathcal{I}_s}),
\]
where ${\rm Sym}(TM)$ denotes the symmetric tensor algebra over $TM$ and $\mathcal{I}_s$ is the ideal spanned by the characteristic polynomial $s$ of $\star$, acting on $E$. Note that in appropriate coordinates, $\star$ is pointwise nothing else than the 'creation' operator of the quantum mechanical harmonic oscillator and the 'diagonalizing' vectors are 'coherent states'.\\

The emphasis of this summarizing discussion will be Hamiltonian systems (cf. \cite{kleinham1}, \cite{kleinham2}), while the examination of the above Lagrangian case and its Frobenius structure, i.e. its connection to 'higher Maslov classes' and miniversal deformations of holomorphic functions with isolated singularities will be continued in (\cite{kleinlag}, see also \cite{kleinham1}). It will turn out that a given Hamiltonian function $H:M\times[0,1]\rightarrow \mathbb{R}$ on a symplectic manifold which is a contangent bundle $(M=T^*N, \omega)$ (we will always assume that the time one map of the corresp. Hamiltonian flow has only non-degenerate fixed points and is of the form $|p|^2$ outside of some compact subset in $T^*N$ containing $N$) also defines a Frobenius structure $\star:TU\rightarrow End(E)$ in analogy to the above, where $E$ is a complex line bundle on a neighbourhood $U$ of the diagonal $\Delta$ in $(M\times M, \omega\oplus\omega)$ so that the corresponding spectral Lagrangian lies in the complexification $(T_{\mathbb{C}}^*U, \omega_{\mathbb{C}})$ and $\pi:L\subset T_{\mathbb{C}}^*U \rightarrow U$ has degree one as well as a $S^1$-valued 'generating function' on $U\subset M\times M$ in the above sense. This function $\Theta$ can be considered to live on $U\subset M \times M$ since $L$ has degree one, then the critical points of $\Theta$ on $U$ correspond exactly to the fixed points of the time-one map of the Hamiltonian flow on $M\times M$, where one extends the Hamiltonian flow of $H$ to $M\times M$ by taking $\tilde H(x,y)=1/2(H(x)+H(y))$ on $U$ (we will assume that $|d\Theta|\rightarrow \infty$ near the boundary of $U$). Since the critical points of the generating function $\Theta$ on $U$ also correspond to the zeros of the spectral Lagrangian, we have the theorem
\begin{theorem}\label{theorem1}
A Hamiltonian function $H:M\times[0,1]\rightarrow \mathbb{R}$ on a cotangent bundle $M=T^*N$ as above defines a Frobenius structure $\star:TU \rightarrow End(E)$ over a neighbourhood $U$ of the diagonal of $(M\times M, \omega \oplus \omega)$, $E$ being a complex line bundle over $U$, so that the following discrete subsets in $U$ coincide:
\begin{itemize}
\item the intersection of the spectral Lagrangian $L$ in $T_{\mathbb{C}}^*U$ with the zero section in $T_{\mathbb{C}}^*U$.
\item the fixed points of the time one flow of $\tilde H$ on $U$. 
\item the critical points of the corresponding generating function $\Theta: U\rightarrow \mathbb{C}^*$.
\end{itemize}
These points are in turn in bijective correspondence to the fixed points of the time one flow of $H$ on $M$.
\end{theorem}
Note that the latter correspondence follows simply by choosing $U$ sufficiently small and altering $\tilde H$ outside $\Delta\subset U$ so that its only fixed points lie on $\Delta$. Note further that we have to pass to a neighbourhood of the diagonal $U\subset M\times M$ to identify the critical points of an $S^1$-valued function $\Theta$ with the fixed points of the time-one flow of $H$ for reasons which will become clear in \cite{kleinham2} (it is closely connected to the question of finding invariant Lagrangian subspaces for the differential of the time one flow of $H$). A Frobenius structure and a spectral Lagrangian living in the complex bundle $T_{\mathbb{C}}^*M$ is always associated to $H$ on $M$ alone, but the zeros of the corresponding spectral Lagrangian do not necessarily correspond to the critical points of a function on $M$ given by matrix elements associated to $E$ (as opposed to the case of the Frobenius structure associated to a 'real' Lagrangian of degree one as above), while they still coincide with the fixed points of the time one flow of $H$. In fact in general and alternatively to the above consideration of a Frobenius structure over $U\subset M\times M$, one has to consider a certain 'dual' Frobenius structure $E'$ to $E$ in the sense that matrix elements associated to $E'$ give (logarithmic) generating functions for the spectral Lagrangian associated to $E$ and vice versa, cf. \cite{kleinham1}. Note also, that for general $M=T^*N$, we have to embed $N$ into a higher dimensional affine space $A$ using the embedding theorem of Nash and Moser (a certain almost complex structure on $TM$ determining the embedding) and then proceed by pulling back the symplectic spinor bundle over $T^*A\times T^*A$ to $U\subset M\times M$. We will give a discussion for $N=T^n$, where $T^n$ denotes the flat torus, which requires no such embedding, then $\Theta$ is again determined by special theta values. Note finally that the spectral Lagrangian $L$ in $T^*_{\mathbb{C}}U$ is not connected to the image of the zero section in $T^*N$ under the time one flow of $H$ in an obvious way.\\
To estimate the number of fixed points of the time one flow of $\tilde H$ on $U$, note that the class $\xi=\Theta^*(\frac{dz}{z}) \in H^1(U, \mathbb{R})$ associated to $\Theta:U\rightarrow \mathbb{C}^*$ defines a local system $\mathcal{L}_{\xi}$ over $U$ by the ring homomorphism
\[
\phi_{\xi}:\mathbb{Z}[\pi]\rightarrow {\bf Nov}(\pi), \quad \phi_{\xi}(g)=t^{<\xi,g>}
\]
where $\pi=\pi_1(U)=\pi_1(M)$ is the fundamental group, $\mathbb{Z}[\pi]$ its group ring, ${\bf Nov}(\pi)$ is the Novikov ring in the indeterminate variable $t$ and $<\xi,g> \in \mathbb{R}$ denotes the evaluation of $\xi$ on the homology class represented by $g \in H_1(U, \mathbb{Z})$. $\mathcal{L}_{\xi}$ is then a left ${\bf Nov}$-module over $U$. Recall that the Novikov ring denotes formal sums 
\[ 
\sum_{i=1}^{\infty}n_it^{\gamma_i},
\]
where $\gamma_i \in \mathbb{R}, \gamma_i\rightarrow -\infty$ and $n_i \in \mathbb{Z}$ are unequal to zero for only a finite number of $i$ obeying $\gamma_i>c$ for any given $c\in \mathbb{R}$. Let $b_i(\xi)$ denote the rank of $H_i(U; \mathcal{L}_{\xi})$ as a module over ${\bf Nov}(\pi)$ and $q_i(\xi)$ the minimal number of generators of its torsion part. Then by the Novikov inequalities resp. their generalizations to manifolds with boundary (cf. Bravermann \cite{bravermann}), Theorem \ref{theorem1} allows to estimate the number of geometrically distinct critical points of $\Theta$ and thus the number of fixed points of $H$ on $M$ by
\begin{folg}
Let $\phi_H$ be the time-one flow of a time-dependent Hamiltonian $H$ on $M$, $n={\rm dim}{M}$ and $\#{\rm Fix}(\phi_H)$ be the number of its fixed points. Then we have the following estimate:
\[
\#{\rm Fix}(\phi_H)\geq \sum_{i=0}^{2n} b_i(\xi)+2 \sum_{i=1}^{2n} q_i(\xi)+ q_0(\xi).
\]
\end{folg}
We assume here that $\Theta$ is modified along a tubular neighbourhood of the boundary $\partial U$ to match the conditions in \cite{bravermann} (which can always be achieved without introducing new critical points). Note that the Novikov numbers $b_i(\xi), q_i(\xi)$ equivalently appear as Betti- resp. torsion numbers of the $\mathbb{Z}[\pi_1(U)]$-module $H_i(\tilde U_{\xi}, \mathbb{Z})$ on the covering $\tilde U_{\xi}$ of $U$ associated to the kernel of the monodromy homomorphism $Per_\xi:\pi_1(U)\rightarrow \mathbb{R}, \ [\gamma]\mapsto <\gamma, \xi>$. Here, $\pi_1(U)$ act as the group deck-transformations on $\tilde U_{\xi}$. We expect to extract further information on the critical points of $\Theta$ by examining the structure of the underlying Morse-Novikov-complex on the chain level more closely. In especially, in the absence of 'homoclinic orbits' estimates involving Lusternik-Schnirelman-like categories of the type introduced in Farber (\cite{farber}) give estimates like the following.
\begin{folg}
Let $\phi_h$ be the time-one flow of a time-dependent Hamiltonian $H$ on $M$ as above and let ${\rm cat}(U, \xi)$ be the category of $U$ with respect to $\xi$ as in introduced in Farber \cite{farber}. Assume that the homology class $[\xi] \in H^1(U, \mathbb{R})$ admits a gradient-like vector field with no homoclinic cycles. Then
\[
\#{\rm Fix}(\phi_H)\geq {\rm cat}(U, \xi).
\]
\end{folg}
Now following the concept of Viterbo \cite{viterbo} and Oh \cite{oh}, we are tempted to define spectral invariants associated to $\Theta$ on $U$ as follows. Denote by $C_*(\tilde U_{\xi})$ the simplicial or cellular chain complex on $\tilde U_{\xi}$, then the Novikov complex $C_*$, generated by the critical points of $\xi$ on $U$ over ${\bf Nov}(\pi)$ is represented as $C_*={\bf Nov}(\pi)\otimes_{\mathbb{Z}[\pi]}C_*(\tilde U_{\xi})$. Let $\Theta_{\xi}:\tilde U_{\xi} \rightarrow \mathbb{R}$ be a primitive of $\xi$ on $\tilde U_{\xi}$. For $\alpha \in C_*$, represent $\alpha =\sum_{i=1}^{\infty}n_{[p,g]} t^{<\xi,g>}$, where $p$ is a critical point of $\Theta$, $g \in \pi$ and $<\xi,g> \in \mathbb{R}$ is the period mapping. We define the level $\lambda_{\xi}(\alpha)$ of $\alpha \in C_*$ as
\[
\lambda_{\xi}(\alpha)=\max_{[p,g]}\{\Theta_{\xi}([p,g]):\ n_{[p,g]}\neq 0\}
\]
Note that $\Theta_{\xi}([p,g])=\Theta(p)+<\xi,g>$ by the definition of the covering $\tilde U_{\xi}$. $\lambda_{\xi}$ defines a filtration on $C_*$ by considering $C_*^\lambda$ as the span of all chains $\alpha$ so that $\lambda_{\xi}(\alpha)\leq\lambda$. There is a natural inclusion $i_\lambda:C_*^\lambda\rightarrow C_*$ and an associated map on $H_*(\tilde U_{\xi}, \mathbb{Z})$. Then we define for any $a \in H_*(\tilde U, \mathbb{Z})$:
\[
\rho(H, a)=\inf_{\alpha; (i_\lambda)[\alpha]=a} \lambda_{\xi}(\alpha).
\]
Note that for $\rho(H, a)$ be finite, necessarily $a \neq 0$, so unless we guarantee the existence of some non-zero homology class $a$ in $H_*(\tilde U_{\xi}, \mathbb{Z})$, we cannot prove the finiteness of $\rho(H, a)$. However, we will prove in this article the following finiteness, spectrality and $C_0$-continuity-property, further investigations and applications of this spectral invariant are postponed to a subsequent paper.
\begin{theorem}
Assume there is a non-zero, non-torsion element in $H_*(\tilde U_{\xi}, \mathbb{Z})$ being a module over ${\bf Nov}(\pi)$. Then $\rho(H, a)$ is finite and a critical value of $\Theta_{\xi}$ for any $0\neq a \in H_*(\tilde U_{\xi}, \mathbb{Z})$. Furthermore, if $H$ and $F$ are two (time-dependent) Hamiltonian fucntions, then
\[
|\rho(H, a)-\rho(F, a)|\leq ||H-F||,
\]
where $||\cdot||$ is Hofer's pseudo-norm on $C_0(T^*N\times [0,1])$. I.e., $\rho_a$ mapping $H\mapsto \rho(H, a)$ is $C_0$-continuous.
\end{theorem}
Note that the construction of such a spectral invariant for a Hamiltonian system on a general cotangent bundle $T^*N$ here goes beyond the reach of Viterbo's finite dimensional methods in \cite{viterbo}, which are in the Hamiltonian case only applicable for $T^*N=\mathbb{R}^{2n}$. The proof of the above finiteness and $C_0$-continuity property leans very closely to the existing proofs of Viterbo and Oh in their respective contexts. This is possible since our 'generating function' $\Theta$ can be interpreted as a 'crude version' of Chaperon's method of broken geodesics resp. Conley and Zehnder's proof of the Arnol'd conjecture for flat tori. However, we want to stress that the main objective of this paper was not to give sharper lower bounds for the existence of Hamiltonian fixed points on cotangent bundles, but to show that the notion of Frobenius stuctures and fundamental questions of symplectic topology are very intimately connected. Interpreting $\Theta$ at least for the case of the torus $M=T^*T^n$ as assuming 'special values' of a certain automorphic function following Mumford's remarks \cite{mumford}, the  connection given in Theorem \ref{theorem1} between the spectral cover of a Frobenius structure associated to the vector bundle $E$ and the critical points of $\Theta$ should have an interpretation in the realms of the Langlands program as giving some sort of 'characteristic zero' analogy for the correspondence between Galois representations and automorphic representations. In especially, the relation between the two complex line bundles $E$ and $\mathcal{L}_{\xi}$ above deserves a closer examination. To both sides, the 'Galois representation side' (the action of the Hamiltonian flow resp. $E$) and the 'automorphic side' (the gradient like-flow of $\Theta$ resp. $\mathcal{L}_{\xi}$) one can associate a dynamical zeta-function (cf. \cite{hutchings}), both should be in a sense 'dual' to another. We mention in particularly recent results of Deligne and Flicker \cite{deligneflicker}. Thus our 'Frobenius structure side' seems to correspond to the Galois-representation side in their article. On the other hand, invoking a dynamical zeta function associated to our $\Theta_{\xi}$ as introduced by Hutchings \cite{hutchings} and associated Lefschetz Theorems one should be able to arrive at a remotely similar theorem as given in their paper, replacing the Frobenius action on the respective objects by the action of the Hamiltonian flow.\\
We finally formulate a conjecture which connects the above spectral invariants (if nontrivial) with the 'eigenvalues' of the covariant derivative of the Euler vector field $X_E$ (which is globally defined on $\tilde U_{\xi}$) associated to the Frobenius structure $\star:TU \rightarrow End(E)$ over $U$ for a non-degenerate Hamiltonian $H$ on $M$. Note that 'eigenvalues' we call here the evaluation of the closed part of the one form with values in $End(E)$ associated to $\nabla X_E$ on a set of generators of $H_1(U, \mathbb{Z})$, in the case of the Frobenius structure associated to the miniversal deformation of an isolated singularity these will be expected to coincide with the usual definition of spectrum (cf. \cite{kleinlag}). The non-triviality of such a closed part follows once one assumes $\xi \in H^1(U, \mathbb{R})$ is non-trivial and $H^*(M, \mathbb{C})$ is {\it formal}, that is all higher order cohomology operations vanish. Note further that our construction of $\star$ should associate a 'variation of Hodge structure' to any Hamiltonian $H$ on a cotangent bundle by the common scheme (cf. \cite{fernan}) of interpreting Frobenius manifolds in terms of 'variations of Hodge structure' and vice versa. On the other hand, our generating function $\Theta$ should be linked to a 'Gromov-Witten'-type theory and its variation of Hodge structures by selecting topologically 'relevant' coherent subbundles of $i^*\mathcal{Q}$ over $M$ by a Thom-isomorphism and thus defining a Frobenius structure on $H^*M$ (cf. a subsequent publication). In any case, we conjecture here, complementing Theorem \ref{theorem1}:
\begin{conj}
The 'eigenvalues' (in the above sense) of $\nabla X_E$ over $\tilde U_{\xi}$, that is the spectrum of the Frobenius structure $\star:TU \rightarrow End(E)$ (that is the spectral numbers of the variation of Hodge structures associated to $H$) coincide generically (after eventual affine scaling) with the above spectral numbers $\rho(H, a)$ of $H$, where $a$ ranges over all elements $a \in H_*(\tilde U_{\xi}, \mathbb{Z})$.
\end{conj}
Note that together with Theorem \ref{theorem1} and interpreting our function $\Theta$ as the kernel of an appropriate integral operator and invoking a related trace formula, this conjecture should be interpreted as an analogon of the (conjectural) Hecke eigenvalue/Frobenius eigenvalue correspondence in the (geometric) Langlands program, an analogous result will be examined in the second paper of this series (\cite{kleinlag}).\\
We want to thank the IHES at Bures sur Yvette, where parts of this research was done, for support and kind hospitality.

\end{document}